\documentclass[twoside, 11pt]{article}
\usepackage{mathrsfs,amsfonts,amsmath}
\RequirePackage{ifthen,calc}
\usepackage{theorem}
\usepackage[dvips]{color}

 \renewcommand{\baselinestretch}{1.0}

 \catcode`@=11
 \def\@evenhead{\hbox to\textwidth{\footnotesize\rm\thepage \hfill
  {\it Yuqiang LI}}} % authors name

 \def\@oddhead{\hbox to \textwidth{\footnotesize{\it
 Riemann-Liouville processes and  Branching Systems } \hfill\thepage}}% abbreviate title

 \renewcommand{\section}{\makeatletter
 \renewcommand{\@seccntformat}[1]{{\csname the##1\endcsname.}\hspace{0.45em}}
 \makeatother \@startsection
{section}%                                            the name
{1}%                                                  the level
{0pt}%                                                the indent
{\baselineskip}%                                      the beforeskip
{0.5\baselineskip}%                                   the afterskip
{\normalsize\bfseries\mathversion{bold}}}
\catcode`@=12

\newcommand\acks{\section*{Acknowledgements}}
\newtheorem{thm}{\noindent Theorem}[section]

\newtheorem{prop}{\noindent Proposition}[section]

\theoremheaderfont{\normalfont\bfseries} \theorembodyfont{\slshape}
{\theorembodyfont{\rmfamily}} {\theorembodyfont{\rmfamily}
}
{\theorembodyfont{\rmfamily} } {\theorembodyfont{\rmfamily}
}
{\theorembodyfont{\rmfamily} } {\theorembodyfont{\rmfamily}
\newtheorem{rem}{\noindent Remark}[section]}
{\theorembodyfont{\rmfamily} } {\theorembodyfont{\rmfamily} } {\theorembodyfont{\rmfamily}
}
 \def\beqlb{\begin{eqnarray}}\def\eeqlb{\end{eqnarray}}
 \def\beqnn{\begin{eqnarray*}}\def\eeqnn{\end{eqnarray*}}
 \numberwithin{equation}{section}
\def\qed{\hfill$\square$\smallskip}
\def\R{{\mathbb R}}
\def\bfE{{\mathbb{E}}}

\def\e{\mathrm{e}}

\def\1int{\int_{0+}^\infty}
\def\0int{\int_0^\infty}
\def\d{\mathrm{d}}

\def\l{\langle}
\def\r{\rangle}
\begin{document}

\title{\LARGE\bf Riemann-Liouville processes arising from Branching particle systems
}
\author{Yuqiang LI \renewcommand{\thefootnote}{\fnsymbol{footnote}}\footnotemark[1]
\\ \small School of Finance and Statistics, East China Normal University,
\\ \small Shanghai 200241, P. R. China.
}
\date{}

\maketitle
\renewcommand{\thefootnote}{\fnsymbol{footnote}}\footnotetext[1]
{Research is supported partially by NSFC grant (10901054) and ECSC grant (13zz037).}

\begin{abstract}

\noindent{It is proved in this paper that Riemann-Liouville processes can arise from the temporal structures of the scaled occupation time fluctuation limits of the site-dependent $(d,\alpha,\sigma(x))$
branching particle systems in the case of $1=d<\alpha<2$ and  $\int_{\R}\sigma(x)\d x<\infty$.}

\smallskip


\smallskip


\end{abstract}

\renewcommand{\baselinestretch}{1.0}

\normalsize

\section{Introduction}
A Riemann-Liouville process (RL process) is a centered Gaussian process. It can be represented by the moving average representation of Brownian motion
 \beqlb\label{Def:RL}
 \int_0^t (t-u)^{H-1/2} B(\d u), \qquad 0\leq t<\infty.
 \eeqlb
Here $B(t)$ denotes the standard Brownian motion and $H>0$ is a constant. For convenience, we denote the RL process by $\{R^H(t),t\geq0\}$. It is easy to see from (\ref{Def:RL}) that $R^H$  does not have stationary increments but is self-similar with index
$H$. RL processes are very useful and important because they are closely  related to
the Riemann-Liouville derivative of Gaussian noise. As early as 1950s, L\'{e}vy \cite{L53} had
briefly commented on this kind of processes. Since then, much work has been done in order to study it; see, for example, \cite{H09,LLS06,SL95}
and the references therein.

From (\ref{Def:RL}) one can readily think of the fractional Brownian motion (FBM) $B^H$ with
(Hurst) index $H\in(0,\:1)$. The FBM $B^H$  has  a moving average representation as follows.
 \beqlb\label{Def:fBm1}
 B^H (t):=\int_{-\infty}^t\Big[ (t-u)_+^{H-1/2}-(-u)_+^{H-1/2}\Big] B(\d u),
 \eeqlb
where $u_+=u\vee 0$. It is well-known that $B^H$ is a
centered Gaussian process with the
covariance function
 \beqlb\label{Def:fBm}
 \bfE[B^{H}(t)B^{H}(s)]=\frac{1}{2}(|s|^{2H}+|t|^{2H}-|t-s|^{2H}).
 \eeqlb
and $B^{H}$ is self-similar with index $H$ and has stationary increments.

Sub-fractional Brownian motions (Sub-FBMs) were formally
introduced in 2004 independently by  Bojdecki et al \cite{BGT04} and Dzhaparidze and van Zanten \cite{DZ04},
respectively. The Sub-FBM $\{\eta^H(t), t\geq0\}$ with index $H\in(0, 1)$ is also a
centered Gaussian process with the covariance function
 \beqlb\label{Def:sub-fBm}
 \bfE[\eta^{H}(t)\eta^{H}(s)]=|s|^{2H}+|t|^{2H}-\frac{1}{2}((s+t)^{2H}+|t-s|^{2H}).
 \eeqlb
This kind of processes has many properties which are the same as those of FBMs, but its increments are not stationary.

When $H\not=1/2$, RL processes, FBMs and Sub-FBMs share the same property that their increments on
non-overlapping intervals are correlated. In fact, Bojdecki et al
\cite{BGT04} showed that the
covariance of increments on intervals at distance $t$ decays like $t^{2H-3}$
for the Sub-FBM $\eta^H$ and  $t^{2H-2}$ for the FBM $B^H$. In addition, as we will point out in Proposition 2.1,  this index for the RL process
$R^H$ is $t^{H-3/2}$.

 Bojdecki et al \cite{BGT04,BGT061} studied the functional limits of the occupation time fluctuation processes
 \beqlb\label{def:oc}
 X_T(t)=\frac{1}{F_T}\int_0^{Tt}\big(N(s)-\bfE(N(s))\big)\d s,\quad t\geq0,
 \eeqlb
of a classical $(d, \alpha,1)$-branching particle system, where $N(t)$ denotes its random counting measure at time $t$, and found that %if the system has
%no branching, then under the condition $d<\alpha$, $F_T=T^{1-\frac{d}{2\alpha}}$ and the
%FBM $B^H$ with $H=1-\frac{d}{2\alpha}$ arises as the temporal structure of the limit process of
%$X_T$ and that
under the assumptions $\alpha<d<2\alpha$ and
$F_T=T^{\frac{3}{2}-\frac{d}{2\alpha}}$, a Sub-FBM $\eta^H$ with $H=\frac{3}{2}-\frac{d}{2\alpha}$ can arise from the temporal structure of the limit process of $X_T(t)$. Here $\bfE(N(s))$ is the
expectation function understood as $\langle
\bfE(N(s)),\phi\rangle=\bfE(\langle N(s), \phi\rangle)$ for any
$\phi\in\mathcal{S}(\R^d)$, the space of smooth rapidly decreasing
functions. Motivated by the aforementioned works, the main purpose of this paper is to report a similar but new result. We show that under certain
conditions, RL processes arise from the limits of occupation time
fluctuations of the so-called site-dependent branching particle systems (see Li \cite{L11a} for details). This result is not obvious at a first sight. In fact, so far as we know, it is not formally
reported in existing literature even if RL processes have been explored for longer time
than that of FBMs and Sub-FBMs.  %The results of this paper provide an interesting complement of the discussion in Bojdecki et al \cite{BGT04,BGT061}.

Precisely, in this paper, we consider a branching particle system in $\R$, whose particles  start off at time $t=0$ from
a Poisson random field with Lebesgue intensity measure $\lambda$, evolve independently with a symmetric $\alpha$-stable L\'{e}vy motion and undergo at rate $\gamma$ the critical, finite variance branching which depends on the position $x$ where the particle splits and is controlled by the moment generating function
 $$g(s,x)=s+\sigma(x)(1-s)^2, \qquad 0\leq s\leq1, 0\leq \sigma(x)\leq 1/2, x\in\R,$$
where $\sigma(\cdot)$ is a measurable function on $\R$. For convenience, we refer this model
to the $(1, \alpha, \sigma(x))$-particle system. Let $N(s)$ denote the random counting measure of the $(1, \alpha, \sigma(x))$-particle system. Consider
the scaled occupation time fluctuation process defined by (\ref{def:oc}). We prove that when $\alpha\in(1, 2)$,
and $F_T=T^{\frac{3}{2}-\frac{1}{\alpha}}$, the RL process $R^H$ with
$H=\frac{3}{2}-\frac{1}{\alpha}$ arises from the temporal structure of the fluctuation limit process if $\int_{\R}\sigma(x)\d x<\infty$. In addition, we also prove that the standard Brownian motion can arise from the limit process  in the case of $\alpha=1$ if $F_T=T^{1/2}\ln T$ and $\int_{\R}\sigma(x)\d x<\infty$.

There is much other literature in the field of functional limit of occupation time fluctuations
of branching systems. For example, Li et al \cite{LX09,L10} studied the branching particle systems with degenerate branching.
Milo\'{s} \cite{M07} investigated the change caused by the varied initial distributions of the classical
$(d, \alpha,\beta)$ branching particle systems.

Without other statement, in this paper, we use $M$ to denote an
unspecified positive finite constant which may not necessarily be
the same in each occurrence, and use $B, B^H, \eta^H, R^H$ to denote the standard Brownian motion,
the FBM with index $H$, the Sub-FBM with index $H$ and the RL process with index $H$, respectively.

The remainder of this paper is organized as follows. In Section 2,
we report the main results. Section  3 contains some basic calculations. In the last
section, i.e. Section 4, we prove the main results.

\section{Main results}

%Let $\vec{\alpha}=(\alpha_1,\alpha_2,\cdots,\alpha_d)$ and denote
%the diagonal matrix
%$diag(\alpha_1^{-1},\alpha_2^{-1},\cdots,\alpha_d^{-1})$ by
%$H(\vec{\alpha})$. By convention,
%$r^{H}=\sum_{k=0}^\infty\frac{(H\ln r)^k}{k!}$ for any $r>0$ and
%$d\times d$ matrix $H$. Therefore
%$$r^{H(\vec{\alpha})}=diag(r^{\alpha_1^{-1}},r^{\alpha_2^{-1}},\cdots,r^{\alpha_d^{-1}}).$$
%Suppose $N=\{N(t), t\geq0\}$ is a $(d, \vec{\alpha},
%\sigma(x))$-branching particle system. The corresponding space
%motion is denoted by $\vec{\xi}$. Then $i-$th component of
%$\{\vec{\xi}(t), t\geq0\}$, $\xi_i(t)$, is a symmetric
%$\alpha_i$-stable Levy process. $\vec{\xi}$ is an
%operator-self-similar process with independent increments (see \cite{S91}).
From now on, we always let $N(s)$ and $X_T(t)$ be the random counting measure of the $(1, \alpha, \sigma(x))$-particle system and the corresponding scaled occupation time fluctuation process. Let $\xi$ denote the particles' space motion. Then $\xi$ is the symmetric $\alpha$-stable L\'{e}vy process. Denote its semigroup by
$\{L_t\}_{t\geq0}$ and the transition density by $p_t$, i.e.,
 $$L_tf(x):=\bfE(f(\xi(t+s))|\xi(s)=x)=\int_{\R}p_t(y-x)f(y)\d y,$$
for all $s, t\geq0$, $x\in\R$ and bounded measurable functions
$f$ (to avoid ambiguity we sometimes write $L_tf(x)$ as
$L_t(f(\cdot))(x)$). For every $\varphi\in\mathcal{S}(\R)$, denote the Green potential operator by
 $$G\varphi(x)=\int_0^\infty L_s\varphi(x)\d s,$$
and define $$G_t\varphi(x)=\int_0^t L_s\varphi(x)\d s.$$
From Li \cite{L11a},
 \beqlb\label{s2-3}
 \bfE(\langle N(t), \phi\rangle)=\int_{\R}L_t\phi(x)\d
 x=\int_{\R} \phi(x)\d x=\langle \lambda, \phi\rangle.
 \eeqlb
The occupation time fluctuation process $X_T=\{X_T(t),
t\geq0\}$ is rewritten as follows.
 \beqlb\label{s2-4}
 \langle X_T(t),\phi\rangle=\frac{1}{F_T}\int_0^{Tt}\langle N(s)-\lambda,
 \phi\rangle\d s,
 \eeqlb
for every $\phi\in\mathcal{S}(\R)$, where $F_T$ is a suitable
norming constant.

%Before we state our main results, we need a concept of convergence; see also \cite{BGT072}.
%
%Let $X, X_T$, $T>0$ be $\mathcal{S}'(\R^d)$-valued c\`{a}dl\`{a}g processes.
%We say that the laws of $X_T$ converge to $X$ in the space-time, or integral, sense
%(denoted $X_T\Rightarrow_i X $) if, for each $t>0$, $\int_0^t\langle X_T, \psi(\cdot, s)\rangle\d s$
%converge to $\int_0^t\langle X, \psi(\cdot, s)\rangle\d s$ in distribution for all $\psi\in\mathcal{S}(\R^{d+1})$.

The main results of this paper read as follows.

\begin{thm}\label{s2-thm-2}
Suppose $D:=\int_{\R}\sigma(x)\d x<\infty$ and
$\alpha\in(1, 2)$. Let
 $F_T^2=T^{3-2/\alpha}$. Then $X_T\Rightarrow K\lambda R^H$ in $C([0, \tau], \R)$ for any $\tau>0$,
where $R^H(\cdot)$ is a RL process with index $H=3/2-1/\alpha$ and $K=\sqrt{2\gamma D}\Gamma(1/\alpha)/[\pi(\alpha-1)]$.
\end{thm}

\begin{thm}\label{s2-thm-1}
 Suppose $D:=\int_{\R}\sigma(x)\d x<\infty$ and $\alpha=1$. Let
 $F_T^2=T (\ln T)^2$. Then $X_T\Rightarrow_i C\lambda B,$ where
$ C=2\sqrt{\gamma D}/\pi$. Here $\Rightarrow_i $ means that for each $t>0$, $\int_0^t\langle X_T, \psi(\cdot, s)\rangle\d s$
converge to $\int_0^t\langle X, \psi(\cdot, s)\rangle\d s$ in distribution for all $\psi\in\mathcal{S}(\R^{2})$.
\end{thm}

The long range dependence of processes can be characterized by means of the so-called {\it
dependence exponent} which was developed in \cite{BGT071}. According to \cite{BGT071}, the
dependence exponent $\kappa$ of a Gaussian process $\{\zeta_t;t\geq0\}$ is
 \beqlb\label{def:de}
 \kappa=\inf_{0\leq u<v<s<t}\sup\Big\{\gamma>0:\; \text{Cov}(\zeta_v-\zeta_u, \zeta_{T+t}-\zeta_{T+s})=o(T^{-\gamma})\;\; \text{as}\;\; T\to\infty\Big\}.\quad
 \eeqlb

The following proposition is of the dependence exponent of RL processes. It seems simple; we write down it because we do not find an appropriate reference.

\begin{prop}\label{s2-prop-3}
The dependence exponent $\kappa$ of the RL process $\{R^H(t), t\geq0\}$ with $H\not=\frac{1}{2}$ is $\frac{3}{2}-H$.
\end{prop}

Below we make some comments on our results.
\begin{rem}\label{s2-rm-2}
%(1) Li \cite{L11a,L11b} had shown that for the $(d, \vec{\alpha}, 1,
%\sigma(x))$ branching system, if $D:=\int_{\R^d}\sigma(x)\d x<\infty$ and $\bar{\alpha}>1$, then $F_T=T^{1/2}$
%and the limit processes $X$ are Gaussian processes with covariances
%\beqnn
% &&{\rm Cov} (\langle X(r),\phi_1\rangle,\langle X(t), \phi_2\rangle)=2(r\wedge t)\int_{\R^d}\Big(\phi_1(x)G\phi_2(x)
% +\gamma\sigma(x)G\phi_1(x)G\phi_2(x)\Big)\d
% x.
% \eeqnn
%Therefore, under the precondition $D:=\int_{\R^d}\sigma(x)\d x<\infty$, only the case of $d=1$
%and $\alpha=2$ is left, on which we will discuss in other place.

(1) From the assumption of $\int_{\R}\sigma(x)\d x<\infty$, one can readily relate the $(1,\alpha,\sigma(x))$-particle systems to the particle systems without branching. Some results of the latter's occupation time fluctuations were revealed in Bojdecki et al \cite{BGT061,BGT062}. They considered the particle systems where particles start off at time $t=0$ from
a Poisson random field with Lebesgue intensity measure $\lambda$ and move independently in $\R^d$ with a symmetrical $\alpha$-stable L\'{e}vy motion, and proved that if $d/\alpha<1$, then $F_T=T^{1-\frac{d}{2\alpha}}$
and the limit process $X=K\lambda B^{1-\frac{d}{2\alpha}}$; if $d/\alpha=1$, then $F_T=(T\ln T)^{1/2}$
and the limit process $X=K\lambda B$. % and if $d/\alpha>1$, then $F_T=T^{1/2}$ and the limit process
%$X$ is a Gaussian process with covariance function
% $${\rm Cov} (\langle X(r),\phi_1\rangle,\langle X(t), \phi_2\rangle)=2(r\wedge t)\int_{\R^d}\phi_1(x)G\phi_2(x) \d x.$$
%(Though they discuss the classical $(d,\alpha,1)$-branching systems, it is easy to see that all
%results in \cite{BGT061,BGT062} hold for $(d,\vec{\alpha},1)$-branching systems with some simple
%modifications.) Naturally, one will think that the particle system $(d, \vec{\alpha}, \sigma(x))$  with $\int_{\R^d}\sigma(x)\d x<\infty$ is similar to the particle system without
%branching.
Comparing their results with the corresponding results in this paper, there at least exist two differences in the case of $d=1$. They are, respectively, $F_T=T^{1/2}\ln T$ vs $F_T=(T\ln T)^{1/2}$ when $\alpha=1$, and  $F_T=T^{\frac{3}{2}-\frac{1}{\alpha}}$ plus $R^{\frac{3}{2}-\frac{1}{\alpha}}$ vs $F_T=T^{1-\frac{1}{2\alpha}}$ plus $B^{1-\frac{1}{2\alpha}}$
 when $\alpha\in(1, 2)$.
The limit processes $R^{\frac{3}{2}-\frac{1}{\alpha}}$ and $B^{1-\frac{1}{2\alpha}}$ have the same dependence
exponent $1/\alpha$ but different self-similarity indice: $\frac{3}{2}-\frac{1}{\alpha}$ vs $1-\frac{1}{2\alpha}$,
and different increments' property: non-stationary increments vs stationary increments. The difference on the norming $F_T$
also shows that the occupation time fluctuations  in the $(1, \alpha, \sigma(x))$ model increase faster than those in the particle systems without branching.

(2) It is also interesting to compare the results of $(1, \alpha, \sigma(x))$-particle systems where $\sigma(x)\equiv 1/2$
with the results in this paper. In fact, when $\sigma(x)\equiv 1/2$, the $(1, \alpha, \sigma(x))$-particle systems are the classical $(d,\alpha,\beta)$-branching particle systems with $d=1,\beta=1$ (see, for example, Bojdecki et al
\cite{BGT061}). %ose of classical $(d, \alpha, 1)$- branching particle systems (cf the both Theorem 2.2 in Bojdecki et al
%\cite{BGT061,BGT062}). As we had mentioned in Li \cite{L11b}, the jumps caused by $\bar{\alpha}=2$
%in $(d, \vec{\alpha}, 1)$ branching systems are smoothed by the condition
%$\int_{\R^d}\sigma(x)\d x<\infty$ in our $(d, \vec{\alpha}, \sigma(x))$ models. For $(d, \vec{\alpha}, 1)$ branching systems, when $\bar{\alpha}\in (1, 2)$
For the classical $(1,\alpha, 1)$-branching particle systems, Theorem 2.2 in Bojdecki et al \cite{BGT061} (see, also, Bojdecki et al \cite{BGT04} ) told us that a long-range dependent Sub-FBM $\eta^{\frac{3}{2}-\frac{1}{2\alpha}}$
with dependence exponent $1/\alpha$ arise from the temporal structure of the limit process of occupation time fluctuations in the case of $\alpha\in(1/2, 1)$. Moreover,  Bojdecki et al intuitively explained  the cause of the long range
dependence as the ``clan recurrence". In this paper, we get the long-range dependent RL
process $R^{\frac{3}{2}-\frac{1}{\alpha}}$ in the case of $\alpha\in(1, 2)$. The long-range dependence exponent is $1/\alpha<1$, and the long-range
dependence is likely caused by the recurrence of particles' motion.

(3) By the terminology in Bojdecki et al \cite{BGT072}, the convergence $\Rightarrow_i $ in Theorem 2.2 is called convergence in the space-time, or integral, sense. Bojdecki et al \cite{BGT86} pointed out that $\Rightarrow_i $ convergence resembles the convergence of finite-dimensional distributions but neither implies the other. However, by the same proof as in Bojdecki et al \cite{BGT072}, we can additionally conclude in Theorem 2.2 that $X_T$ converges to $C\lambda B$ in the sense of finite-dimensional distributions.

\end{rem}

%In the sequel, we will make use of the following lemma, see Remark 2.1 in \cite{LX09}.
%\begin{lem}\label{s2-lm-1}
%Let $z=(z_1,\cdots, z_d)$. For $r\in(0, \bar{\alpha})$,
% $ \int_{[0,
%1]^d}\frac{1}{\sum_{k=1}^d|z_k|^{r\alpha_k}}\d z<\infty,
% $
%and for $r>\bar{\alpha}$,
% $ \int_{\R^d\setminus[-1,
%1]^d}\frac{1}{\sum_{k=1}^d|z_k|^{r\alpha_k}}\d z<\infty.
% $
% Furthermore, if $\phi(z)$
%is bounded and $\int_{\R^d} \phi(z)\d z<\infty$, then
%  $ \int_{\R^d}\frac{\phi(z)}{\sum_{k=1}^d|z_k|^{r\alpha_k}}\d
% z<\infty,$
%for all $r\in(0, \bar{\alpha})$.
%\end{lem}

At the end of this section, let us show some heuristics to explain why one can expect the RL process arising as a limit in Theorem \ref{s2-thm-2}. To this end, observe the covariance of $X_T$. Let $0\leq s, t\leq 1$. From (\ref{tight0}), it is easy to get that
 \beqnn
&&Cov(\l X_T(t),\phi\r, \l X_T(s),\phi\r)=\frac{T^2}{F_T^2}\int_0^s\int_0^t Cov\Big(\l N(Tu),\phi\r, \l N(Tv),\phi\r\Big)\d u \d v\nonumber
\\&&\qquad\qquad=\frac{2\gamma T^3}{F_T^2} \int_0^s\int_0^t\int_{\R}\sigma(x)\int_0^{u\wedge v}L_{T(u-r)}\phi(x) L_{T(v-r)}\phi(x)\d r\d x\d u \d v\nonumber
\\&&\qquad\qquad\quad+\frac{T^2}{F_T^2}\int_0^s\int_0^t\int_{\R}\phi(x) L_{T|u-v|}\phi(x)\d x\d u\d v.
\eeqnn
From  (3.1) and (3.2) in Bojdecki et al \cite{BGT071}, for any $\phi\in\mathcal{S}(\R)$, there is  constant $M>M_1>0$ such that
 \beqlb\label{tight4}
 M_1 u^{-1/\alpha}\leq L_u\phi(x)\leq Mu^{-1/\alpha}
 \eeqlb
for any $u>0$. Since  $F_T^2=T^{3-\frac{2}{\alpha}}$ and $1>1/\alpha$, we have that
 \beqnn
 &&\frac{\gamma T^3}{F_T^2} \int_0^s\int_0^t\int_{\R}\sigma(x)\int_0^{u\wedge v}L_{T(u-r)}\phi(x) L_{T(v-r)}\phi(x)\d r\d x\d u \d v\nonumber
 \\&&\qquad\qquad \sim \int_0^s\int_0^t\d u\d v\int_0^{u\wedge v} (u-r)^{-1/\alpha}(v-r)^{-1/\alpha}\d r\nonumber
 \\&&\qquad\qquad \sim \int_0^{s\wedge t} (s-r)^{1-1/\alpha}(t-r)^{1-1/\alpha}\d r,
 \eeqnn
where $A\sim B$ means that $A/B\in(0,\infty)$, and that
  $$\lim_{T\to\infty}\frac{T^2}{F_T^2}\int_0^s\int_0^t\int_{\R}\phi(x) L_{T|u-v|}\phi(x)\d x\d u\d v=0.$$
Therefore,
 $$Cov(\l X_T(t),\phi\r, \l X_T(s),\phi\r)\sim  \int_0^{s\wedge t} (s-r)^{1-1/\alpha}(t-r)^{1-1/\alpha}\d r,$$
which is the covariance of the RL process $R^H$ with $H=3/2-1/\alpha$.

%In the rest of this paper,
% $\widehat{f}(z)$, $z\in\R^d$, denotes the Fourier transform of integrable function $f$, i.e.,
% $$\widehat{f}(z)=\int_{\R}\e^{ixz}f(x)\d x.$$
%$L_t$, as the semi-group of a symmetric $\alpha$-stable L\'{e}vy process, satisfies that
% \beqlb\label{s2-7}
% \widehat{L_t f}(z)=\widehat{f}(z)\e^{-t|z|^{\alpha}}.
% \eeqlb
%It is well-known that if $\phi\in\mathcal{S}(\R)$, then $\widehat{\phi}(z)$ is bounded and integrable and goes to $0$ as $|z|\to\infty$.

%We remind that if $\phi\in \mathcal{S}(\R)$ and $\alpha\in(1, 2)$, then by the inverse Fourier transform formula
%    \beqnn
%  \frac{|G_T\phi(x)|}{T^{1-\bar{\alpha}}}&=&\frac{1}{T^{1-\bar{\alpha}}}\Big|\int_0^T L_u\phi(x)\d u\Big|
%  \leq\frac{1}{(2\pi)^d}\frac{1}{T^{1-\bar{\alpha}}}\int_{\R^d}\frac{|\widehat{\phi}(z)|(1-\e^{-T \sum_{i=1}^d|z_i|^{\alpha_i}})}{\sum_{i=1}^d|z_i|^{\alpha_i}}\d z
%  \\&\leq&\frac{1}{(2\pi)^d}\int_{\R^d}\big|\phi(x)\big|\d x\int_{\R^d}\frac{1-\e^{- \sum_{i=1}^d|z_i|^{\alpha_i}}}{\sum_{i=1}^d|z_i|^{\alpha_i}}\d z<\infty.
%  \eeqnn
%  Hence, there is a constant $K>0$ such that for any $x\in\R^d$ and $T>0$,
%   \beqlb\label{s2-9}
%   |G_T\phi(x)|<KT^{1-\bar{\alpha}}.
%   \eeqlb

\section{Preliminary calculations}
Define a sequence of random variables
$\tilde{X}_T$ in $\mathcal{S}'(\R^2)$ as follows: For any
$T>1$ and $\psi\in\mathcal{S}(\R^2)$, let
 \beqlb\label{s3-0}
 \langle\tilde{X}_T, \psi\rangle=\int_0^1\langle X_T(t), \psi(\cdot, t)\rangle\d t.
 \eeqlb
Without other statement, in the sequel, $\psi\in\mathcal{S}(\R^2)$
always has the form $\psi(x, t)=\phi(x)h(t)$, where
$\phi\in\mathcal{S}(\R)$ and $h\in\mathcal{S}(\R)$ are
nonnegative functions. Let
 \beqlb\label{s3-2}
 \tilde{h}(s)=\int_s^1 h(t)\d t\;\;\text{and}\;\;\;\psi_T(x,
s)=\frac{1}{F_T}\phi(x)\tilde{h}(\frac{s}{T}).
 \eeqlb
Define
 \beqlb\label{s3-1-2}
 V_{\psi_T}(x, t, r)=1-\bfE_x\Big(\exp\Big\{-\int_0^t\langle N(s), \psi_T(\cdot, r+s)\rangle\d s\Big\}\Big),
 \eeqlb
where $\bfE_x(f(N)):=\bfE(f(N)|N(0)=\epsilon_x)$ and $\epsilon_x$ denotes the unit measure concentrated at $x\in\R$. From Section 3 in \cite{L11a}, we get that
 \beqlb\label{s3-1-5'}
 V_{\psi_T}(x, t, r)&=&\int_0^tL_s\big[\psi_T(\cdot, r+s)\big(1-V_{\psi_T}(\cdot, t-s, r+s)\big)\big](x)\d s\nonumber
\\ &&-\gamma \int_0^tL_s[\sigma(\cdot)V_{\psi_T}^2(\cdot, t-s,
 r+s)](x)\d s.
 \eeqlb
 Furthermore, define
  \beqlb\label{s3-1-4}
 J_{\psi_T}(x, t,
 r):=\int_0^t L_s\psi_T(\cdot, r+s)(x)\d s.
 \eeqlb
Then
 \beqlb\label{s3-1-4'}
 V_{\psi_T}(x, t,
r)\leq J_{\psi_T}(x, t,
 r),\quad
\eeqlb
and
  \beqlb\label{s3-1-6}
   \bfE(\e^{-\langle \tilde{X}_T, \psi\rangle})=\exp\Big(I_1(T,\psi_T)+I_2(T,\psi_T)+I_3(T,\psi_T)\Big),
   \eeqlb
 where
   \beqlb
   I_1(T,\psi_T)&=&\gamma\int_{\R}\sigma(x)\d x\int_0^T J_{\psi_T}^2(x, T-s, s)\d s,\label{s3-1-10}
   \\I_2(T,\psi_T)&=&\gamma\int_{\R}\sigma(x)\d x\int_0^T [V_{\psi_T}^2(x, T-s, s)-J_{\psi_T}^2(x, T-s, s)]\d s\nonumber
   \\&\geq& -2\gamma
 \Big[I_{21}(T,\psi_T)+\gamma I_{22}(T,\psi_T)\Big],\qquad\label{s3-1-11}
   \\ I_3(T,\psi_T)&=&\int_{\R}\d x\int_0^T\psi_T(x,
  s)V_{\psi_T}(x, T-s, s)\d s,\label{s3-1-8}
   \eeqlb
and
 \beqlb
 I_{21}(T,\psi_T)&=&\int_{\R}\sigma(x)\d x\int_0^TJ_{\psi_T}(x, T-s, s)\d
 s\nonumber
 \\&&\times\int_s^TL_{u-s}\big(\psi_T(\cdot, u)J_{\psi_T}(\cdot, T-u,
 u)\big)(x)\d u,\qquad\label{s3-3-2}
\\  I_{22}(T,\psi_T)&=&\int_{\R}\sigma(x)\d x\int_0^TJ_{\psi_T}(x, T-s, s)\d
 s\nonumber
 \\&&\times\int_s^T L_{u-s}(\sigma(\cdot) J_{\psi_T}^2(\cdot, T-u, u))(x)\d
 u.\qquad\label{s3-3-3}
 \eeqlb

\medskip

\noindent{\bf 3.1 Limits of $I_1$, $I_2$ and $I_3$ when  $1<\alpha<2$}

\medskip

{\bf Step 1}. We consider the limit of $I_1(T,\psi_T)$. From (\ref{s3-1-10}) and (\ref{s3-1-4}), we get that
 \beqlb\label{s3-B-8}
 I_1(T,\psi_T)&=&\gamma\int_{\R}\sigma(x)\d x\int_0^T\bigg(\int_0^{T-s}L_u\psi_T(x, s+u)\d
 u\bigg)^2\d
 s\nonumber
\\&=&\int_0^1h(t)\d t\int_0^1 h(r)\d r\int_0^{r\wedge t}\Psi_T(s, r, t)\d s,
 \eeqlb
where
  \beqlb\label{s3-B-9}
  \Psi_T(s, r, t)=\frac{T\gamma}{F_T^2}\int_{\R}\sigma(x)\d x\int_0^{T(t-s)}L_u\phi(x)\d u\int_0^{T(r-s)}L_v\phi(x)\d v.
  \eeqlb
The self-similarity of the $\alpha$-stable L\'{e}vy process $\xi$ implies that
 $$\frac{1}{T^{1-\frac{1}{\alpha}}}\int_0^{Tt}L_u f(x)\d u=\int_0^t\int_{\R} u^{-1/\alpha}p_1((x-y)T^{-1/\alpha}u^{-1/\alpha})f(y)\d y,$$
for any integrable function $f$. This means that for every $0<t\leq 1$ and non-negative integrable function $f$
 \beqlb\label{s3-B-4}
\frac{\alpha}{\alpha-1}t^{1-\frac{1}{\alpha}}p_1(0)\int_{\R} f(y)\d y\geq \frac{1}{T^{1-\frac{1}{\alpha}}}\int_0^{Tt}L_u f(x)\d x.
 \eeqlb
Additionally, as $T\to\infty$,
\beqnn
\frac{1}{T^{1-\frac{1}{\alpha}}}\int_0^{Tt}L_u f(x)\d x\to\frac{\alpha}{\alpha-1}t^{1-\frac{1}{\alpha}}p_1(0)\int_{\R}f(y)\d y.
\eeqnn
Therefore, (\ref{s3-B-9}) and the facts that $\int_{\R}\sigma(x)\d x<\infty$ and $F_T^2=T^{3-2/\alpha}$ lead to
  %\beqnn
%  \Psi_T(s)&=&\frac{\gamma}{(2\pi)^{2}T^{2-2/\alpha}}\int_{\R^{2}}\overline{\widehat{\sigma}(z+w)}
%  \widehat{\phi}(z)\widehat{\phi}(w)\nonumber
%  \\&&\times\Big[\frac{1-\e^{-{T(t-s)|z|^{\alpha}}}}{|z|^{\alpha}}
%  \frac{1-\e^{-{T(t-s)|w|^{\alpha}}}}{|w|^{\alpha}}\Big]\d z\d w.
%  \eeqnn
%Then by using the substitutions $z'=T^{1/\alpha}z$ and $w'=T^{1/\alpha}w$, we get that
 \beqnn
  \Psi_T(s, r, t)
  &\to&\gamma (p_1(0)\frac{\alpha}{\alpha-1})^2[(t-s)(r-s)]^{1-1/\alpha}\int_{\R}\sigma(x)\d x\Big(\int_{\R}\phi(x)\d x\Big)^2
  \\&=& C[(t-s)(r-s)]^{1-1/\alpha}\Big(\int_{\R}\phi(x)\d x\Big)^2,
  \eeqnn
as $T\to\infty$, where
 \beqnn
 C&=&\gamma \Big(\frac{\alpha p_1(0)}{\alpha-1}\Big)^2\int_{\R}\sigma(x)\d x=\gamma \Big(\frac{\Gamma(1/\alpha)}{\pi(\alpha-1)}\Big)^2\int_{\R}\sigma(x)\d x.
  \eeqnn
Consequently, from (\ref{s3-B-8}) we can readily have that as $T\to\infty$,
 \beqlb\label{s3-B-10}
  I_1(T,\psi_T)\to C\Big(\int_{\R}\phi(x)\d x\Big)^2\int_0^1h(t)\d t\int_0^1 h(r)
  \d r\int_0^{r\wedge t}[(t-s)(r-s)]^{1-1/\alpha}\d s.\qquad
 \eeqlb

{\bf Step 2}. We consider the limit of $I_2(T,\psi_T)$. Using (\ref{s3-2}), (\ref{s3-1-4}) and the fact that $\tilde{h}$ is bounded, we have that
 \beqlb\label{s3-B-1}
 J_{\psi_T}(x, T-u, u)\leq \frac{M}{F_T}G_T\phi(x),
 \eeqlb
for any $u\leq T$ and some $M>0$. Applying
this inequality to (\ref{s3-3-2}) and using (\ref{s3-2}) again,
% \beqnn\label{s3-B-2}
% I_{21}(T,\psi_T)&=&\frac{1}{F_T^3}\int_{\R}\sigma(x)\int_0^T\d s
%\int_0^{T-s}\tilde{h}\Big(\frac{s+u}{T}\Big)\d u
% \int_0^{T-s}\tilde{h}\Big(\frac{s+v}{T}\Big)\d
% v\nonumber
% \\&&\times\int_0^{T-s-u}\bigg[\tilde{h}\Big(\frac{s+u+w}{T}\Big)L_v\phi(x)L_u(\phi L_w\phi)(x)
%\bigg]\d
% x.
% \eeqnn
%Due to $h\in\mathcal{S}(\R)$, .
we get a constant $M>0$ such that
 \beqlb\label{s3-B-3}
  I_{21}(T,\psi_T)&\leq&\frac{M T}{F_T^3}
\int_{\R}\sigma(x)\int_0^TL_v\phi(x)\d v
 \int_0^TL_u(\phi\int_0^TL_w\phi\d
w)(x)\d u\d x\nonumber
\\&\leq&\frac{M T}{F_T^3}\int_{\R}\sigma(x)G_T\phi(x)G_T(\phi G_T\phi)(x)\d x.
\eeqlb
It follows from (\ref{s3-B-4}) that for any non-negative integrable function $f$ there exists a constant $M>0$ such that
\beqlb\label{s2-9}
G_T f(x)<M T^{1-1/\alpha}.
\eeqlb
Therefore, substituting $F_T^2=T^{3-2/\alpha}$ into (\ref{s3-B-3}), we can readily obtain that
 \beqlb
    I_{21}(T,\psi_T)\leq \frac{M}{T^{1/2}}\int_{\R}\sigma(x)\d x\to 0.
 \eeqlb
Furthermore, applying (\ref{s3-B-1}) and (\ref{s3-2}) to (\ref{s3-3-3}), we
obtain that
% \beqnn
%I_{22}(T,\psi_T)&=&\int_{\R}\sigma(x)\int_0^T\d
%s\int_0^{T-s}\d u\int_0^{T-s-u}\d t\int_0^{T-s-u}\d
%t'\int_0^{T-s}\d v
%\\&&\times\Big[L_t\psi_T(x, s+u+t)L_{t'}\psi_T(x,
%s+u+t')L_u(\sigma(\cdot)L_v\psi_T(\cdot, s+v))(x)\Big]\d x,\nonumber
% \eeqnn
%where we use the fact $\int_{\R}fL_ug\d x=\int_{\R}gL_uf\d x$ for all bounded integral functions $f, g$.
%Then by the same discussion as above there exists a constant $M>0$ such that
  \beqlb\label{s3-B-5}
I_{22}(T,\psi_T)&\leq&\frac{M T}{F_T^3}%\int_{\R}\sigma(x)\int_0^T\d u\int_0^{T-u}
%L_t\phi(x)\d t\int_0^{T-u} L_{t'}\phi(x)\d t'\nonumber
%\\&&\qquad\qquad\qquad\times\Big[ L_u\Big(\sigma\int_0^T L_v\phi\d v\Big)(x)\Big]\d x.\qquad
% \eeqlb
%By (\ref{s2-9}) again, there is a constant $M>0$,
% \beqlb\label{s3-B-5-1}
% &&\int_{\R}\sigma(x)\int_0^TL_t\phi(x)\d t\int_0^T L_{t'}\phi(x)\d t'\int_0^T L_u\Big(\sigma\int_0^T L_v\phi\d v\Big)(x)\d u\d
%x\nonumber
%\\&&\quad\leq
\int_{\R}\sigma(x)(G_T\phi(x))^2G_T(\sigma G_T\phi)(x)\d x,%\leq M T^{3(1-1/\alpha)}\int_{\R}\sigma(x)G_T\sigma(x)\d x.
 \eeqlb
for some $M>0$. From (\ref{s2-9}), (\ref{s3-B-5}) and the fact $\int_{\R}\sigma(x)\d x<\infty$, it follows that
% $$ I_{22}(T,\psi_T)\leq\frac{M}{T^{1/2}}\int_{\R}\sigma(x)G_T\sigma(x)\d x\leq\frac{M}{T^{1/2}}
% \int_{\R}\frac{1-\e^{-T|z|^\alpha}}{|z|^\alpha}\d z.$$
%Let $y=T^{1/\alpha}z$. Then for $\alpha\in (1, 2)$,
 \beqlb\label{s3-B-6}
 I_{22}(T,\psi_T)\leq\frac{M}{T^{1/\alpha-1/2}}\to 0,
 \eeqlb
where the convergence is due to the fact $\alpha\in (1, 2)$. Note that $I_2(T,\psi_T)\leq 0$ which follows from (\ref{s3-1-4'}). Combining (\ref{s3-1-11}) with (\ref{s3-B-3}) and (\ref{s3-B-6}) yields that as
$T\to\infty$,
 \beqlb\label{s3-B-7}
 I_2(T,\psi_T)\to 0.
 \eeqlb

{\bf Step 3}. At last we pass to the limit of $I_3(T,\psi_T)$. Let
  \beqlb\label{s3-B-11}
 I_{31}(T, \psi_T):=\int_{\R}\d x\int_0^T\psi_T(x,
 s)J_{\psi_T}(x, T-s, s)\d s.
 \eeqlb
 From (\ref{s3-1-8}), (\ref{s3-1-4'}), (\ref{s3-B-1}) and (\ref{s2-9}), it follows that
 \beqlb\label{s3-B-12}
 0\leq I_3(T,\psi_T)&\leq& I_{31}(T,\psi_T)%=\int_{\R}\Big[\int_0^T\d s\int_0^{T-s}\psi_T(x,
% s)L_v\psi_T(x, s+v)\d v\Big]\d x.\qquad
% \eeqlb
%By (\ref{s3-2}) and (\ref{s2-9}),
%  \beqlb\label{s3-B-13}
% I_3(T,\psi_T)
%&\leq&\frac{M}{F_T^2}\int_{\R}\phi(x)\int_0^T\d s\int_0^{T-s}
% L_v\phi(x)\d v\d
% x\nonumber
 \leq\frac{M T}{F_T^2}\int_{\R}\phi(x)G_T\phi(x)\d
 x\nonumber
 \\&\leq&\frac{T^{2-1/\alpha}M}{F_T^2}\int_{\R}\phi(x)\d x,
 \eeqlb
for some $M>0$. Substituting $F_T^2=T^{3-2/\alpha}$ into (\ref{s3-B-12}), $\alpha\in (1, 2)$ indicates
 that
  \beqlb\label{s3-B-14}
  I_3(T,\psi_T)\to 0.
  \eeqlb

\noindent{\bf 3.2 Limits of $I_1$, $I_2$ and $I_3$ when $\alpha=1$}

 \medskip

%In this case, $d=1,\vec{\alpha}=\alpha_1=1$ or $d=2, \vec{\alpha}=(\alpha_1,\alpha_2)=(2,2)$.

We first point out that
 \beqlb\label{s3-A-1}
 \lim_{T\to\infty}I_2(T,\psi_T)=0,
 \eeqlb
and
 \beqlb\label{s3-A-2}
 \lim_{T\to\infty}I_3(T,\psi_T)=0.
 \eeqlb
The details are similar to those of $I_2$ and $I_3$
in Section 3.1, and hence we omit them.  Below, we discuss $I_1(T,\psi_T)$.

%In fact, from (\ref{s3-1-10}) it follows that
% \beqlb\label{s3-A-8}
% I_1(T,\psi_T)&=&\gamma\int_{\R}\sigma(x)\d x\int_0^T\bigg(\int_0^{T-s}L_u\psi_T(x, s+u)\d
% u\bigg)^2\d
% s\nonumber
%\\&=&\int_0^1h(t)\d t\int_0^1 h(r)\d r\int_0^{r\wedge t}\Psi_T(s, r,t)\d s,
% \eeqlb
%where
%  \beqnn
%  \Psi_T(s, r,t)=\frac{T\gamma}{F_T^2}\int_{\R}\sigma(x)\d x\int_0^{T(t-s)}L_u\phi(x)
%  \d u\int_0^{T(r-s)}L_v\phi(x)\d v.
%  \eeqnn
From (3.47) in \cite{BGT082} we know that for any $0<t\leq 1$,
 \beqlb
 \lim_{T\to\infty}\frac{1}{\ln T}\int_0^{Tt}L_u\phi(x)\d u=p_1(0)\int_{\R}\phi(x)\d x.
 \eeqlb
Using this fact and noting that $\int_{\R}\sigma(x)\d x<\infty$ and $F_T^2=T(\ln T)^2$,  we get  that
  \beqlb\label{s3-A-9}
  \lim_{T\to\infty}\Psi_T(s,r,t)&=&\gamma \Big[p_1(0)\int_{\R}\phi(x)\d x\Big]^2\int_{\R}\sigma(x)\d x,
  \eeqlb
where $\Psi_T(s,r,t)$ is the same as (\ref{s3-B-9}).
Note that (\ref{s2-9}) and the integrability of $\sigma(x)$ imply the boundedness of $\Psi_T(s, r,t)$ with respect to $T>0$, $0\leq r, t\leq 1$ and $0\leq s\leq r\wedge t$.    Therefore, the dominated convergence theorem and (\ref{s3-B-8}) imply that as $T\to\infty$,
 \beqlb\label{s3-A-11}
  I_1(T,\psi_T)\to C_1\Big(\int_{\R}\phi(x)\d x\Big)^2\int_0^1h(t)\d t\int_0^1 h(r) (r\wedge t)\d r,\qquad
 \eeqlb
where
 $$ C_1=\gamma(p_1(0))^2\int_{\R}\sigma(x)\d x
 =\frac{\gamma }{\pi^{2}}\int_{\R}\sigma(x)\d x.$$

\section{Proofs of the main results}

\medskip

\noindent{{\bf Proof of Theorem 2.1} Without loss of generality, fix $\tau=1$. As explained in Bojdecki et al \cite{BGT061},%By the space-time method in Bojdecki et al. \cite{BGT86},
to prove this theorem it is sufficient to show that
 \begin{itemize}
 \item[(i)] For each non-negative $\psi\in\mathcal{S}(\R^{d+1})$, as $T\to\infty$
 \beqlb\label{s4-1}
 \bfE(\e^{-\langle\tilde{X}_T, \psi
 \rangle})\to\exp\bigg(\frac{1}{2}\int_0^1\int_0^1 {\rm Cov}\Big(\big\langle X(s), \psi(\cdot, s)\big\rangle,\big\langle X(t),\psi(\cdot, t)\big\rangle\Big)\d s\d
 t\bigg),
 \eeqlb
 where $\tilde{X}_n$ and $\tilde{X}$ are defined as
 (\ref{s3-0}) and $X$ is the corresponding limit process;
 \item[(ii)] $\{\langle X_T, \phi\rangle; T\geq 1\}$ is tight in $C([0, 1],
 \R)$ for any given $\phi\in\mathcal{S}(\R)$, where the
theorem of Mitoma \cite{M83} is used.
\end{itemize}
If $\psi(x, t)=\phi(x)h(t)$, where
$\phi\in\mathcal{S}(\R)$ and $h\in\mathcal{S}(\R)$ are
nonnegative functions, then the limit process in Theorem 2.1 satisfies that
 \beqnn
 &&{\rm Cov}\Big(\big\langle X(s), \psi(\cdot, s)\big\rangle,\big\langle X(t),\psi(\cdot, t)\big\rangle\Big)=K^2\langle \lambda, \phi\rangle^2 h(s)h(t)\bfE(R(s)R(t)),
 \\&&\qquad\qquad=K^2\langle \lambda, \phi\rangle^2 h(s)h(t)\int_0^{s\wedge t}[(t-u)(s-u)]^{1-1/\alpha}\d u,
 \eeqnn
where
 $K^2=2\gamma \int_{\R}\sigma(x)\d x \big(\frac{\Gamma(1/\alpha)}{\pi(\alpha-1)}\big)^2.$
Furthermore, combining (\ref{s3-1-6}) with (\ref{s3-B-7}), (\ref{s3-B-10}) and (\ref{s3-B-14})
yields
 \beqnn
 &&\lim_{T\to\infty}\bfE(\e^{-\langle\tilde{X}_T, \psi
 \rangle})
 \\&&=\exp\bigg(\gamma \Big(\frac{\Gamma(\frac{1}{\alpha})\langle \lambda, \phi\rangle}{\pi(\alpha-1)}\Big)^2\int_{\R}\sigma(x)\d x
 \int_0^1\int_0^1 h(s)h(t)\int_0^{s\wedge t}[(t-u)(s-u)]^{1-\frac{1}{\alpha}}\d u\d s\d t\bigg)\nonumber
  \\&&=\exp\bigg(\frac{1}{2}\int_0^1\int_0^1 {\rm Cov}\Big(\big\langle X(s), \psi(\cdot, s)\big\rangle,\big\langle X(t),\psi(\cdot, t)\big\rangle\Big)\d s\d
 t\bigg).
 \eeqnn
This means that (\ref{s4-1}) holds for the special case $\psi(x, t)=\phi(x)h(t)$.

For general $\psi$, the proof of (i) is the same
with slightly more complicated notation. The  details are omitted.

 We now pass  to prove (ii). Note that by some simple and standard arguments, one can readily get that for any $f, g\in \mathcal{S}(\R)$ and $s\leq t$,
 \beqlb\label{tight0}
 \bfE_x\Big(\l N(s),f\r \l N(t),g\r\Big)=L_s( f L_{t-s}g)(x)+2\gamma\int_0^s L_u(\sigma L_{s-u}f L_{t-u}g)(x)\d u.
 \eeqlb
Thus, by Poisson initial condition,
$$Cov\Big(\l N(s),f\r, \l N(t),g\r\Big)=\int_{\R}f(x) L_{t-s}g(x)\d x+2\gamma\int_{\R} \sigma(x) \int_0^sL_{s-u}f(x) L_{t-u}g(x)\d u\d x.$$
From the setting of $X_T(t)$, we have that for any $0\leq s<t\leq 1$,
\beqlb\label{tight1}
&&\bfE(\l X_T(t),\phi\r-\l X_T(s),\phi\r)^2=\frac{T^2}{F_T^2}\int_s^t\int_s^t Cov\Big(\l N(Tu),\phi\r, \l N(Tv),\phi\r\Big)\d u \d v\nonumber
\\&&\qquad\qquad=\frac{2\gamma T^3}{F_T^2} \int_s^t\int_s^t\int_{\R}\sigma(x)\int_0^{u\wedge v}L_{T(u-r)}\phi(x) L_{T(v-r)}\phi(x)\d r\d x\d u \d v\nonumber
\\&&\qquad\qquad\quad+\frac{T^2}{F_T^2}\int_s^t\int_s^t\int_{\R}\phi(x) L_{T|u-v|}\phi(x)\d x\d u\d v.
\eeqlb
Applying the last inequality of (\ref{tight4}) to the last two formulas of (\ref{tight1}), noting that $F_T^2=T^{3-\frac{2}{\alpha}}$ and $\sigma(x)$ is integrable, we get that for all $T>1$
 \beqlb\label{tight2}
   &&\frac{2\gamma T^3}{F_T^2} \int_s^t\int_s^t\int_{\R}\sigma(x)\int_0^{u\wedge v}L_{T(u-r)}\phi(x) L_{T(v-r)}\phi(x)\d r\d x\d u \d v\nonumber
   \\&&\qquad\leq M\int_s^t\d u\int_s^t\d v\int_0^{u\wedge v}(u-r)^{-1/\alpha}(v-r)^{-1/\alpha}\d r\nonumber
   \\&&\qquad\leq M\int_s^t\d u\int_s^t|u-v|^{-1/\alpha}\d v\nonumber
   \\&&\qquad=M(t-s)^{2-\frac{1}{\alpha}}\int_0^1\int_0^1|u-v|^{-1/\alpha}\d u\d v\leq M(t-s)^{2-\frac{1}{\alpha}},
 \eeqlb
and
  \beqlb\label{tight3}
 &&\frac{T^2}{F_T^2}\int_s^t\int_s^t\int_{\R}\phi(x) L_{T|u-v|}\phi(x)\d x\d u\d v
 \nonumber
 \\&&\qquad\leq\frac{MT^{2-\frac{1}{\alpha}}}{F_T^2}\int_s^t\int_s^t|u-v|^{-1/\alpha}\d u\d v
 \leq M(t-s)^{2-\frac{1}{\alpha}},
 \eeqlb
hold for constants $M>0$ which only depend on $\phi$.

Combining (\ref{tight1}) with (\ref{tight2}) and (\ref{tight3}), we get that
 $$\bfE(\l X_T(t),\phi\r-\l X_T(s),\phi\r)^2\leq M(\phi)(t-s)^{2-\frac{1}{\alpha}}.$$
Therefore, from Billingsley \cite[Theorem 12.3]{B68} and the fact $\l X_T(0),\phi\r=0$, it follows that $\{\langle X_T, \phi\rangle\}_{T\geq 1}$ is tight in $C([0, 1],  \R)$. The proof of Theorem 2.1 is complete. \qed

\medskip

 \noindent{\bf Proof of Theorem 2.2}. The proof of Theorem 2.2 is easy. In fact, combining (\ref{s3-1-6}) with (\ref{s3-A-1}), (\ref{s3-A-2}) and (\ref{s3-A-11}) we can readily get the corresponding formula (\ref{s4-1}). More details are omitted.\qed

 At last, we give the proof of Proposition 2.1.

 \medskip

\noindent{\bf Proof of Proposition 2.1}. From the representation (\ref{Def:RL}) of $R^H(\cdot)$, it is easy to see that %is
% $$C(r, t)=\int_0^{r\wedge t}[(t-s)(r-s)]^{H-\frac{1}{2}}\d s,\qquad r,t\geq 0.$$
%Suppose $r<t$. Then by the transform $u=s/t$, we get that
% $$C(r, t)=t^{2H}\int_0^{r/t}[(1-u)(\frac{r}{t}-u)]^{H-\frac{1}{2}}\d u=t^{2H}\int_0^{r/t}[(1-\frac{r}{t}+s)s]^{H-\frac{1}{2}}\d s.$$
%For $x\in[0,1]$, let
% $$F(x)=\int_0^x[s(1-x+s)]^{H-\frac{1}{2}}\d s,$$
%and $\Phi(x)=F(x)-\frac{2}{2H+1}x^{H+\frac{1}{2}}.$
%Then by L'H\^{o}pital's rule it is easy to verify that
%  \beqlb\label{prop-2}
% \lim_{x\to0}\frac{\Phi(x)}{x^{H+\frac{3}{2}}}=\lim_{x\to0} \frac{F(x)-\frac{2}{2H+1}x^{H+\frac{1}{2}}}{x^{H+\frac{3}{2}}}=\frac{2(2H-1)}{(2H+3)(2H+1)}.
% \eeqlb Therefore,
for any $0\leq u<v<s<t$,
 \beqlb\label{prop-3}
 &&\text{Cov}\big(R(v)-R(u), R(T+t)-R(T+s)\big)\nonumber
 \\&&=\int_0^{\infty}[(v-r)^{H-\frac{1}{2}}{\bf 1}_{\{r\leq v\}}-(u-r)^{H-\frac{1}{2}}{\bf 1}_{\{r\leq u\}}][(T+t-r)^{H-\frac{1}{2}}-(T+s-r)^{H-\frac{1}{2}}]\d r.\nonumber
 \eeqlb
By Taylor expansion formula, it is easy to see that for any given $t,s, r$,
 $$\Psi(T, t, s,r):=\frac{[(T+t-r)^{H-\frac{1}{2}}-(T+s-r)^{H-\frac{1}{2}}]}{T^{H-\frac{3}{2}}}\to (H-\frac{1}{2})(t-s),$$
as $T\to\infty$. Consequently, one can readily have that
 \beqlb
 &&\lim_{T\to\infty}\frac{\text{Cov}\big(R(v)-R(u), R(T+t)-R(T+s)\big)}{T^{H-\frac{3}{2}}}\nonumber
  \\&&=(H-\frac{1}{2})(t-s)\int_0^{\infty}[(v-r)^{H-\frac{1}{2}}{\bf 1}_{\{r\leq v\}}-(u-r)^{H-\frac{1}{2}}{\bf 1}_{\{r\leq u\}}]\d r\nonumber
  \\&&=\frac{2H-1}{2H+1}\Big(v^{H+\frac{1}{2}}-u^{H+\frac{1}{2}}\Big)(t-s),
 \eeqlb
%By (\ref{prop-2}), we have that for any $v, t\geq 0$,
% \beqlb
% \lim_{T\to\infty}\Phi(\frac{v}{T+t})T^{H+\frac{3}{2}}=\frac{2(2H-1)}{(2H+3)(2H+1)}v^{H+\frac{3}{2}}.
% \eeqlb
%Therefore
% \beqlb
% &&\lim_{T\to\infty}\frac{\text{Cov}\big(R(v)-R(u), R(T+t)-R(T+s)\big)}{T^{H-\frac{3}{2}}}\nonumber
% \\&&=\lim_{T\to\infty}(1+\frac{t}{T})^{2H}T^{H+\frac{3}{2}}\Big(\Phi(\frac{v}{T+t})-\Phi(\frac{u}{T+t})\Big)\nonumber
% \\&&\qquad\qquad-\lim_{T\to\infty}(1+\frac{s}{T})^{2H}T^{H+\frac{3}{2}}\Big(\Phi(\frac{v}{T+s})-\Phi(\frac{u}{T+s})\Big)\nonumber
% \\&&\qquad\qquad+\frac{2}{2H+1}\Big(v^{H+\frac{1}{2}}-u^{H+\frac{1}{2}}\Big)\lim_{T\to\infty}\frac{(T+t)^{H-\frac{1}{2}}-(T+s)^{H-\frac{1}{2}}}{T^{H-\frac{3}{2}}}\nonumber
% \\&&=\
% \eeqlb
which and (\ref{def:de}) together imply that the dependence exponent of RL process $R^H(\cdot)$ with $H\not=\frac{1}{2}$ is $\frac{3}{2}-H$.\qed

\acks
The author thanks the anonymous referees for careful reading of the paper and for the useful
suggestions which have helped to improve the paper significantly.

\end{document}